\documentclass[11pt,openany ]{article}
\usepackage{graphicx}
\usepackage{subfig}
\usepackage{amssymb}
\usepackage{amsmath}
\usepackage{amsthm}
\usepackage{amsfonts}
\usepackage{mathrsfs}
\usepackage{makeidx}
\usepackage{multicol}

\input{cyracc.def}

\numberwithin{equation}{section}

\begin{document}
\title{\Large \bf   On entire solutions for a class of product-type nonlinear PDEs in $\mathbb{C}^n$}
\author{Feng L\"{u}\footnote{Corresponding author, Email: lvfeng18@gmail.com(F. L\"{u})}\\
\small{College of Science, China University of Petroleum, Qingdao, Shandong, 266580, P.R. China. }\\
}
\date{}
\maketitle

\vspace{3mm}

\begin{abstract}
This paper is mainly devoted to describing the entire solutions of nonlinear partial differential equation
$$
u_{z_1}u_{z_2}\cdots u_{z_n}=e^g,
$$
with the eikonal equation as a prototype, where $g$ is a polynomial in $\mathbb{C}^n$. Through a novel method, we break through the restriction of finite order condition and present the explicit expressions for the entire solution of the above equation. As an application, we completely resolve two questions of Xu-Liu-Xuan in \cite{Xu}.
\end{abstract}



\pagenumbering{arabic}

\section{Introduction and main results}
The eikonal equation $\sum_{i=1}^n(\frac{\partial u}{\partial x_i})^2=1$, ($u:\mathbb{R}^n\rightarrow \mathbb{R}$) is a fundamental non-linear partial differential equation (PDE) that arises in the high-frequency asymptotic analysis of wave propagation. Derived from the wave equation (or Maxwell's equations in electromagnetism) using the WKB approximation (geometric optics limit), it serves as a critical bridge between physical wave optics and geometric ray optics, see e.g., \cite{Calin, Courant, Garabedian}. Numerous results have been established for the eikonal equation in the field of real numbers, whereas relatively few studies have been conducted in the field of complex numbers. It seems to us that Khavinson \cite{Khavinson} was the first to focus on the complex-analytic solutions to this equation and investigated the entire solution of the eikonal equation $u_{z_1}^2+u_{z_2}^2=1$. With the linear transformation $x = z_1 + iz_2$ and $y =z_1 - iz_2$, Khavinson transformed the above equation to
\begin{equation}\label{A002}
u_xu_y=1,
\end{equation}
and showed that the entire solution must be a linear function. The result was later derived by Li \cite{Li4} and Saleeby \cite{Saleeby} with different ways, respectively. In 1997, Hemmati applied the method of characteristics to generalize Khavinson's result as follows: Let $u$ be an entire solution in $\mathbb{C}^2$ of $F(u_{z_1},u_{z_2})=0$, where $F$ is an entire function whose zero set $\{(p,q)\in\mathbb{C}^2:F(p,q)=0\}$ does not contain any complex lines, i. e., $F$ does not have a linear factor. Then $u(z_1,z_2)$ is a linear function. Meanwhile, Hemmati also showed that the theorem fails in higher dimensions by the following example, which can also be found in \cite{Khavinson}. Consider $u(z_1,z_2,z_3)=z_3+f(z_1+z_2)$ which solves $u_{z_1}^2-u_{z_2}^2+u_{z_3}-1=0$, where $f$ is an entire function in $\mathbb{C}$. The subsequent contributions are due to Li \cite{Li1,Li2,Li3, Li4}, who described entire solutions of some variations of the eikonal equation with $n=2$, especially this form $u_{z_1}^2+u_{z_2}^2=e^g$ with a polynomial in $\mathbb{C}^2$. Some other studies on variations of the eikonal equation of (\ref{A002}) can be found in \cite{Chen, Lu,Saleeby1}.\\

It can be seen that the equation (\ref{A002}) can be transformed into a Monge-Amp\'{e}re equation by a variable change and differentiation, which is one of the most fundamental and influential fully nonlinear partial differential equations in mathematics. Inspired by this equation (\ref{A002}), Xu-Liu-Xuan  \cite{Xu} considered entire solutions of a more general form of (\ref{A002}) as
\begin{equation}\label{A02}
\prod_{i=1}^3\left(a_{i 1} u_{z_1}+a_{i 2} u_{z_2}+a_{i 3} u_{z_3}\right)=1.
\end{equation}
Let
$$
\mathcal{A}=\begin{pmatrix}
a_{11} & a_{12} & a_{13} \\
a_{21} & a_{22} & a_{23} \\
a_{31} & a_{32} & a_{33}
\end{pmatrix},
$$
where $a_{ij}$ ($i,j=1,2,3$) is a constant. Their result can be stated as follows. \\

\noindent \textbf{Theorem A.} Let $|\mathcal{A}| \neq 0$. Then equation (\ref{A02}) does not have any finite order transcendental entire solution.\\

Here, $|\mathcal{A}|$ is the determinant of $\mathcal{A}$ and the order of a meromorphic function $u$ in $\mathbb{C}^n$ is defined as
$$
\rho(u)=\limsup_{r\rightarrow\infty}\frac{\log T(r,u)}{\log r},
$$
where $T(r,u)$ is the characteristic function in Nevanlinna theory, see e.g., \cite{St, Vi}.\\

In \cite{Xu}, Xu-Liu-Xuan presented several examples to demonstrate that the condition $|\mathcal{A}| \neq 0$ is indispensable. In the same paper, they posed the following two questions.\\

\noindent \textbf{Question 1.} Can the restrictive condition ``finite order" in Theorem A be removed?\\

\noindent \textbf{Question 2.} Is the conclusion of Theorem A still valid for equation (\ref{A02}) with $n(\geq 4)$ complex variables.\\

We are interested in the aforementioned problems. Assume that
\begin{equation}\label{1.1}
\prod_{i=1}^n\left(a_{i 1} u_{z_1}+a_{i 2} u_{z_2}+\cdots+a_{i n} u_{z_n}\right)=1,
\end{equation}
and a matrix
$$
\mathcal{A}=\begin{pmatrix}
a_{11} & a_{12} & \cdots & a_{1n} \\
a_{21} & a_{22} & \cdots & a_{2n} \\
\cdots & \cdots & \cdots & \cdots \\
a_{n1} & a_{n2} & \cdots & a_{nn}
\end{pmatrix}.
$$

Using elementary row transformations of a
matrix, the present author and Ma \cite{Ma} obtained the following result, thereby providing a positive answer to \textbf{Question} 2 under the finite order condition.\\

\noindent \textbf{Theorem B.} Let $|\mathcal{A}| \neq 0$. Then equation (\ref{1.1}) does not have any finite order transcendental entire solution.\\

It should be noted that the finite order condition plays a crucial role in the proof of Theorems A and B. Moreover, the proof procedure in Theorems A and B can hardly be applied to the case of infinite order. In this paper, by adopting a novel method, we overcome the limitation of finite order condition and give a positive answer to the aforementioned problems. More precisely, we derive the following result.\\

\noindent \textbf{Theorem 1.} Suppose that $|\mathcal{A}| \neq 0$ and $u$ is an entire solution of (\ref{1.1}). Then $u$ is a linear function.\\

Prior to the proof of Theorem 1, we investigate the entire solutions of the following equation
 \begin{equation}\label{2.1}
u_{z_1}u_{z_2}\cdots u_{z_n}=e^g,
 \end{equation}
where $g$ is a polynomial in $\mathbb{C}^n$. In order to state the main result, we introduce two definitions.\\

\textbf{Definition 1.} For simplicity, we write $u_{ij}=u_{z_iz_j}$ for any $i,j\in I=\{1,2,...,n\}$. Set
$$
J=\{i\in I: u_{ij}\equiv 0~ \hbox{for any} ~~j\in I\backslash\{i\}\},~~\chi=I\backslash J.
$$
It is pointed out that the set $J$ may be empty. Denote by $\sharp G$ the number of elements in a set $G$.\\

\textbf{Definition 2.} For two equations $u_{z_i}$ and $u_{z_j}$, we define the notations $\sim$ and $\approx$ as follows. If $u_{ij}\not\equiv 0$, we say that $u_{z_i} \sim u_{z_j}$. It follows from $u_{ij}=u_{ji}$ that $u_{z_j} \sim u_{z_i}$ if $u_{z_i} \sim u_{z_j}$. In addition, we say that $u_{z_i} \approx u_{z_j}$ if
$u_{z_i} \sim u_{z_j}$ or $u_{z_i} \sim u_{z_k}\cdots \sim u_{z_j}$. Otherwise, we say that $u_{z_i}\not \approx u_{z_j}$. For example, if $u_{ik}\not\equiv 0,~u_{kt}\not\equiv 0,~u_{tj}\not\equiv 0$, then $u_{z_i}\sim u_{z_k}\sim u_{z_t}\sim u_{z_j}$, which also implies that $u_{z_i}\approx u_{z_j}$.\\

With the above definitions, we derive the following main result. \\

\noindent \textbf{Theorem 2.} Suppose that $u=u(z)=u(z_1,...,z_n)$ is an entire solution of (\ref{2.1}). Then, one of the following assertions holds. \\

(a) If $\sharp J=n$, then $g(z)=\sum_{i=1}^n g_i({z_i})+2 k_0 \pi i$ and
$$
u(z)=\sum_{i=1}^n \int e^{g_i({z_{i}})}dz_{i};
$$
where $k_0$ is a fixed integer.\\

(b) If $\sharp J=0$, then there exists a partition $I=\cup_{l=1}^t I_l$ with $I_l=\{k_{l1},..., k_{l\theta_l}\}$ and $\sharp I_l=m_l\geq 2$, for any $i\in I_\mu$ and $j\in I_\nu$ such that
\begin{equation}\label{0001}
\left\{
\begin{aligned}
	&u_i\approx u_j& \hbox{if}~~1\leq \mu=\nu\leq t, \\
	&u_i\not\approx  u_j ~& \hbox{if}~~1\leq \mu\neq \nu\leq t.\\
\end{aligned}
\right.
\end{equation}
In addition,
$$
g(z)=\sum_{l=1}^t g_l(A_{l1}z_{k_{l1}}+...+A_{l\theta_l}z_{k_{l\theta_l}})+2k_0\pi i,~~u(z)=\sum_{l=1}^t f_l(A_{l1}z_{k_{l1}}+...+A_{l\theta_l}z_{k_{l\theta_l}}),
$$
where $f_l$ is an entire function satisfying $f_l'=e^{\frac{g_l}{m_l}}$, $A_{lj}$ is a constant, and $k_0$ is a fixed integer;\\

(c) If $0<\sharp J=\{\tau_1,...,\tau_\nu\}<n$, then there exists a partition $\chi=I\backslash J=\cup_{l=1}^t I_l$ as in (b) satisfying condition $(\ref{0001})$ (For simplicity, we still use the above notations) such that
$$
\begin{aligned}
&g(z)=\sum_{i=1}^\nu g_i({z_{\tau_{i}}})+\sum_{l=1}^t h_l(A_{l1}z_{k_{l1}}+...+A_{l\theta_l}z_{k_{l\theta_l}})+2k_0\pi i,\\
&u(z)=\sum_{i=1}^\nu \int e^{g_i({z_{\tau_{i}}})}dz_{\tau_{i}}+ \sum_{l=1}^t f_l(A_{l1}z_{k_{l1}}+...+A_{l\theta_l}z_{k_{l\theta_l}}),\\
\end{aligned}
$$
where $f_l$ is an entire function satisfying $f_l'=e^{\frac{h_l}{m_l}}$, $A_{lj}$ is a constant, and $k_0$ is a fixed integer.\\

For the special case $t=1$ in (b) and (c), we have the following corollary.\\

\noindent \textbf{Corollary 1.} Suppose that $u$ is an entire solution of (\ref{2.1}) and $t$ is defined as in (b) and (c) of Theorem 2. Then, we have the following assertions. \\

(1) If $\sharp J=0$ and $t=1$, then
$$
g(z)=g(A_1z_1+A_2z_2+...+A_nz_n),~~u(z)=f(A_1z_1+A_2z_2+...+A_nz_n),
$$
where $f$ is an entire function satisfying that $f'=e^{\frac{g}{n}}$ and $A_1\cdots A_n=1$.\\

(2) If $0<\sharp J=\{\tau_1,...,\tau_\nu\}<n$, $\chi=I\backslash J=\{\tau_{\nu+1},...,\tau_n\}$ and $t=1$, then
$$
\begin{aligned}
&g(z)=\sum_{i=1}^\nu g_i({z_{\tau_{i}}})+g_{\nu+1}(A_{\tau_{\nu+1}}z_{\tau_{\nu+1}}+...+A_{\tau_{n}}z_{\tau_{n}})+2k_0\pi i,\\
&u(z)=\sum_{i=1}^\nu \int e^{g_i({z_{\tau_{i}}})}dz_{\tau_{i}}+ f(A_{\tau_{\nu+1}}z_{\tau_{\nu+1}}+...+A_{\tau_{n}}z_{\tau_{n}}),\\
\end{aligned}
$$
where $f$ is an entire function satisfying $f'=e^{\frac{g_{\nu+1}}{n-\nu}}$, $A_{\tau_j}$ $(j=\nu+1,...,n)$ is a constant, and $k_0$ is a fixed integer.\\

\textbf{Remark 1.} For $n=2,3$, by the above results, we list the forms of entire function $u$ to equation (\ref{2.1}) in the following table.

\begin{table}[h]
\centering
\renewcommand{\arraystretch}{1.6}  
\renewcommand{\tabcolsep}{3pt}  
\footnotesize
\begin{tabular}{|c|c|c|c|}
\hline
$n=2$ & $\sharp J=0$ & $u(z)=f(A_1z_1+A_2z_2)$  &$g(z)=g(A_1z_1+A_2z_2)$\\
\cline{2-4}
      & $\sharp J=2$ & $u(z)=\int e^{g_1(z_1)}dz_1+\int e^{g_2(z_2)}dz_2$ &$g(z)=g_1(z_1)+g_2(z_2)+2k_0\pi i$ \\
\hline
$n=3$ & $\sharp J=0$ & $u(z)=f(A_1z_1+A_2z_2+A_3z_3)$& $g(z)=g(A_1z_1+A_2z_2+A_3z_3)$\\
\cline{2-4}
      & $\sharp J=1$ & $u(z)=\int e^{g_{\tau_1}(z_{\tau_1})}dz_{\tau_1}+ f(A_2z_{\tau_2}+A_3z_{\tau_3})$& $g(z)=g_{\tau_1}(z_{\tau_1})+h(A_2z_{\tau_2}+A_3z_{\tau_3})+2k_0\pi i$\\
\cline{2-4}
      & $\sharp J=3$ & $u(z)=\int e^{g_1(z_1)}dz_1+\int e^{g_2(z_2)}dz_2+\int e^{g_3(z_3)}dz_3$ &$g(z)=g_1(z_1)+g_2(z_2)+g_3(z_3)+2k_0\pi i$\\
\hline
\end{tabular}
\caption{}
\end{table}

\textbf{Remark 2.} In some sense, Theorem 2 can be seen a generalization of \cite[Theorem 2.1]{Li3} given by Li. The function $g$ in Theorem 2 is assumed to be a polynomial. A natural problem is whether or not $g$ can be generalized to be a transcendental entire function in Theorem 2. However, the answer is negative, as shown by the following example. (The method in this paper is not applicable to this problem, which we leave for future research.) \\

\textbf{Example 1.} Consider $u=ie^{\frac{1}{2}(2z_2i+e^{2z_1i})}$ and $g=2(z_1+z_2)i+e^{2z_1i}$. A calculation yields that $u_{z_1}u_{z_2}=e^g$. Obviously, $u$ does not satisfy the conclusion (a)-(c). This example can be seen in \cite[Proposition 2.5]{Li3}.\\

In Theorem 2, it is easy to derive that if $g$ is a constant, then $f_l$ reduces to a linear function. So, the following corollary follows immediately from Theorem 2.\\

\noindent \textbf{Corollary 2.} Suppose that $g$ is a constant in (\ref{2.1}) and $u$ is an entire solution of (\ref{2.1}). Then, $u$ is a linear function.\\

Now, with Corollary 2, we prove Theorem 1.\\

\noindent\textbf{Corollary 2$\Rightarrow$ Theorem 1.} Suppose that $u$ is an entire solution of (\ref{1.1}). Set
\begin{equation}\label{0002}
\begin{pmatrix}
z_{1}  \\
z_{2} \\
\cdots \\
z_{n}
\end{pmatrix}=  \mathcal{A}^T\begin{pmatrix}
w_{1}  \\
w_{2} \\
\cdots \\
w_{n}
\end{pmatrix}.
\end{equation}
Denote $v(w_1,...,w_n)=u(z_1,...,z_n)$. Then, it follows from (\ref{1.1}) that
\begin{equation}\label{2.10}
v_{w_1}v_{w_2}\cdots v_{w_n}=1.
\end{equation}
By Corollary 2, we get that $v(w_1,...,w_n)$ is a linear function in $w_1,...,w_n$. Together with $|\mathcal{A}|\neq 0$, we derive that $u(z_1,...,z_n)$ is a linear function in $z_1,...,z_n$, which is the desired result of Theorem 1.\\

\textbf{Remark 3.} Through the above change of variables $(\ref{0002})$, we can investigate entire solution of
\begin{equation}\label{A003}
\prod_{i=1}^n\left(a_{i 1} u_{z_1}+a_{i 2} u_{z_2}+\cdots+a_{i n} u_{z_n}\right)=e^g,
\end{equation}
where $g$ is a polynomial in $\mathbb{C}^n$ and $|\mathcal{A}| \neq 0$. Here, we omit the details. Furthermore, this method in the present paper may also be applied to analyze the entire solution to equation
 \begin{equation}\label{A00003}
\prod_{i=1}^n\left(a_{i 1} u_{z_1}+a_{i 2} u_{z_2}+\cdots+a_{i n} u_{z_n}\right)=pe^g,
\end{equation}
where $p$ is a polynomial. Naturally, a more refined analysis is required, and the forms of the solutions are more complicated.\\

Before to proceed, we assume the reader's familiarity with the basic notations of Nevanlinna theory and utilize some results in Nevanlinna theory (see e.g. \cite{Hu, St, Vi}). Let $f$ be a meromorphic function in $\mathbb{C}^n$ and $S(r,f)$ denote any quantity satisfying $S(r,f)=o\{T(r,f)\}$ as $r\rightarrow\infty$ outside a set of $r$ of finite Lebesgue measure.

(i). The Nevanlinna first fundamental theorem $T(r,f)=T(r,\frac{1}{f})+O(1)$;

(ii). $T(r,f+g)\leq T(r,f)+T(r,g)+O(1)$, $T(r,fg)\leq T(r,f)+T(r,g)+O(1)$ for any two meromorphic functions $f,~g$ in $\mathbb{C}^n$;

(iii). $T(r, \alpha)=S(r, e^\beta)$, where $\beta$ is a non-constant entire function in $\mathbb{C}^n$ and $\alpha$ is any differential polynomial in $\beta$, including $\beta$, $\beta_{z_j}$ for any $j=1,...,n$.

\section{Proof of Theorem 2}

Suppose that $u$ is an entire solution of (\ref{2.1}). Then
$$
u_{z_1}u_{z_2}\cdots u_{z_n}=e^g.
$$
Observe that $u_{z_i}$ has no zeros and poles. The Hadamard factorization theorem yields that
\begin{equation}\label{2.2}
u_{z_i}=e^{\alpha_i}, ~~(i=1,..,n)
\end{equation}
where $\alpha_i$ is an entire function in $\mathbb{C}^n$ and
\begin{equation}\label{2.3}
e^{\alpha_1+\alpha_2+\cdots+\alpha_n}=e^g,
\end{equation}
which implies that
\begin{equation}\label{2.4}
\alpha_1+\alpha_2+\cdots+\alpha_n=g-2k_0\pi i.
\end{equation}
with a fixed integer $k_0$. Below, we consider two cases.\\

\textbf{Case 1.} $\sharp J=0$. Then $\chi=I$ and the following condition holds:\\

\textbf{Condition}: For any $i$, there exists at least an integer $j\in \{1,2,...,n\}\backslash \{i\}$ such that $u_{ij}\not \equiv0$. \\

We have the following property.\\

\textbf{Property.} If $u_i\approx u_k$, then $e^{\alpha_i}=\gamma_{ik} e^{\alpha_k}$, where $\gamma_{ik}$ is an entire function with $T(r, \gamma_{ik})=S(r,e^{\alpha_k})$. \\

It is emphasis that the property implies $T(r,e^{\alpha_i})=T(r,e^{\alpha_k})+S(r,e^{\alpha_k})$ if $u_i\approx u_k$.\\

Now, we prove the above property. \\

Note that $u_i\approx u_k$. If $u_i=e^{\alpha_i}\sim u_k=e^{\alpha_k}$, one has
$$
\frac{\partial\alpha_i}{\partial z_k}e^{\alpha_i}=u_{ik}=u_{ki}=\frac{\partial\alpha_k}{\partial z_i}e^{\alpha_k}\not\equiv0.
$$
Rewrite it as
$$
e^{\alpha_i}=\frac{\frac{\partial\alpha_k}{\partial z_i}}{\frac{\partial\alpha_i}{\partial z_k}} e^{\alpha_k}=\gamma_{ik} e^{\alpha_k},
$$
which plus the facts (ii) and (iii) leads to the property. Further, if $u_i\sim u_j\sim \cdots \sim u_k$, the same argument also yields the property. So, the property is proved.\\

For the set $I=\{1,...,n\}$, there must exist a partition $I=\cup_{l=1}^t I_l$ for any $i\in I_\mu$ and $j\in I_\nu$ such that
$$
\left\{
\begin{aligned}
	&u_i\approx u_j& \hbox{if}~~1\leq \mu=\nu\leq t, \\
	&u_i\not\approx  u_j ~& \hbox{if}~~1\leq \mu\neq \nu\leq t.\\
\end{aligned}
\right.
$$
Suppose that $\sharp I_l =p_l$ $(1\leq l\leq t)$. The \textbf{Condition} implies that $p_l\geq 2$ $(1\leq l\leq t)$. By rearranging the order $i\in I$ and the corresponding variables $z_i$, without loss of generality, we can assume that
$$
I_1=\{1,...,p_1\}, I_2=\{p_1+1,...,p_1+p_2\},...., I_t=\{p_1+...+p_{t-1}+1,...,p_1+...+p_t=n\}.
$$
For convenience, assume $q_0=0$, $q_1=p_1$ and $p_1+p_2+...+p_\nu=q_\nu$ for $2\leq \nu\leq t$. Then,
$$
I_1=\{1,...,p_1\}=\{q_0+1,...,q_1\}, I_2=\{q_1+1,...,q_2\},...., I_t=\{q_{t-1}+1,...,q_t=n\}.
$$

Below, we consider two subcases.\\

\textbf{Subcase 1.1.} $t=1$. That is $u_i=e^{\alpha_i}\approx u_j=e^{\alpha_j}$ for any $1\leq i\neq j\leq n$.\\

From \textbf{Condition}, for $i=1$, there exists an integer $j\in \{2,...,n\}$ such that $u_{1j}\not\equiv 0$. Without loss of generality, we assume that $u_{12}\not \equiv0$. Then, the fact $u_{12}=u_{21}$ yields that
\begin{equation}\label{2.6}
 \frac{\partial \alpha_1}{\partial z_2}e^{\alpha_1}= \frac{\partial \alpha_2}{\partial z_1}e^{\alpha_2}\not \equiv0.
\end{equation}
Observe that $I_1=I$ and \textbf{Property}. Then, $e^{\alpha_j}=\gamma_{j1} e^{\alpha_1}$ for any $j\in I$, where $\gamma_{j1}$ is an entire function with $T(r,\gamma_{j1})=T(r,e^{\alpha_1})$. Obviously, $\gamma_{11}=1$. Together with (\ref{2.6}), we have
\begin{equation}\label{2.7}
\begin{aligned}
&e^{n\alpha_1}\frac{\partial \alpha_1}{\partial z_2} \prod_{j=3}^{n} \gamma_{j1}
=e^{\alpha_1}\frac{\partial \alpha_1}{\partial z_2}e^{\alpha_1}\prod_{j=3}^{n}\gamma_{j1}e^{\alpha_1} \\
&=e^{\alpha_1}\frac{\partial \alpha_2}{\partial z_1}e^{\alpha_2}\prod_{j=3}^{n}e^{\alpha_j} =\frac{\partial \alpha_2}{\partial z_1}\prod _{j=1}^ne^{\alpha_j} \\
&=e^{\alpha_1+...+\alpha_n} \frac{\partial \alpha_2}{\partial z_1}=e^g \frac{\partial \alpha_2}{\partial z_1},
\end{aligned}
\end{equation}
which implies that
\begin{equation}\label{2.71}
e^{n\alpha_1} \phi=e^g\varphi
\end{equation}
with $\phi=\frac{\partial \alpha_1}{\partial z_2} \prod_{j=3}^{n} \gamma_{j1}$ and $\varphi=\frac{\partial \alpha_2}{\partial z_1}$. Observe that $T(r, \varphi)=S(r, e^{\alpha_1})$, $T(r,\phi)=S(r, e^{\alpha_1})$, which together with the polynomial $g$ leads to that $\alpha_1$ is a polynomial, so are $\alpha_{j}$ $(j=2,...,n)$, $\varphi$ and $\phi$. Furthermore, the equation (\ref{2.71}) yields that $g-n\alpha_1$ is a constant. For any $k\in I$, the same argument yields that $g-n\alpha_k$ is a constant. Set $\alpha_k=\frac{g}{n}+B_k$ with a constant $B_k$ and set $e^{B_k}=A_k$. Then, $u_{z_k}=e^{\alpha_k}=e^{\frac{g}{n}+B_k}=A_ke^{\frac{g}{n}}$ and
$$\frac{u_{z_1}}{A_1}=\frac{u_{z_2}}{A_2}=\cdots=\frac{u_{z_n}}{A_n}=e^{\frac{g}{n}}.$$
A routing way to solve the above PDE yields that
$$
u=f(A_1z_1+A_2z_2+...+A_nz_n),
$$
where $f$ is an entire function satisfying $f'=e^{\frac{g}{n}}$. So, $g=g(A_1z_1+A_2z_2+...+A_nz_n)$. Meanwhile,
$$A_1\cdots A_n=e^{B_1+...+B_n}=e^{\alpha_1+...+\alpha_n-g}=1,$$ which is the conclusion (b).\\

\textbf{Subcase 1.2.} $t\geq 2$.\\

Set $E$ as
$$
E=(u_{ij})=\begin{pmatrix}
u_{1 1} &  ... & u_{1 p_1}& u_{1 (p_1+1)}&...&...&...&u_{1n} \\
u_{2 1} & ... & u_{2 p_1} & u_{2 (p_1+1)}&...&...&...&u_{2n}\\
... &  ... & ... & ...\\
u_{p_1 1}  & ... & u_{p_1 p_1}& u_{p_1 (p_1+1)}&...&...&...&u_{p_1n}\\
u_{(p_1+1) 1}&...&u_{(p_1+1) p_1}&u_{(p_1+1) (p_1+1)}&...&u_{(p_1+1) q_2}&...&...\\
... &  ... & ... & ...&...&...&...&...\\
... &  ... & ... & u_{q_2 (p_1+1)}&...&u_{q_2 q_2}&...&...\\
... &  ... & ... & ...&...&...&...&...\\
u_{n 1} &  ... & ...& ...&...&...&...&u_{nn} \\
\end{pmatrix}.
$$
By the partition $I=\cup_{l=1}^t I_l$, we obtain the fact $u_{\kappa\lambda}=u_{\lambda \kappa}\equiv 0$ for any $1\leq \kappa\leq p_1$ and $p_1+1\leq\lambda\leq n$. Otherwise, $u_\kappa\sim u_\lambda$ and $u_\kappa\approx u_\lambda$, which contradicts $I_1=\{1,...,p_1\}$. This fact implies that
$$
u_{i}=e^{\alpha_i}=e^{\alpha_i(z_1,...,z_{p_1})},~~i=1,...,p_1.
$$
Similarly to the above, we obtain
$$
u_{i}=e^{\alpha_i}=e^{\alpha_i(z_{q_{l-1}+1},...,z_{q_l})},~~i=q_{l-1}+1,...,q_l,~~2\leq l\leq t.
$$
In fact, $E$ is a block diagonal matrix as follows
\[
E=
\begin{pmatrix}
U_1 &  &  &  & \\
& U_2 &  &  & \\
&  & \ddots &  & \\
&  &  & U_l & \\
&  &  &  & \ddots & \\
&  &  &  &  & U_t
\end{pmatrix},
~~U_l=
\begin{pmatrix}
u_{(q_{l-1}+1)(q_{l-1}+1})  & \cdots & u_{(q_{l-1}+1)q_l}\\
u_{(q_{l-1}+2)(q_{l-1}+1)}  & \cdots & u_{(q_{l-1}+2)q_l}\\
\vdots &  \ddots & \vdots\\
u_{q_l(q_{l-1}+1)}  & \cdots & u_{q_lq_l}
\end{pmatrix}, (1\leq l\leq t).
\]

It is easy to deduce that (For simplicity, below we still use the notation $g_i$)
\begin{equation}\label{2.8}
\left\{
\begin{aligned}
	&\alpha_1(z_1,...,z_{p_1})+...+\alpha_{p_1}(z_1,...,z_{p_1})=g_1(z_1,...,z_{p_1}), \\
	&\sum_{i=q_{l-1}+1}^{q_l}\alpha_i=\sum_{i=q_{l-1}+1}^{q_l}\alpha_i(z_{q_{l-1}+1},...,z_{q_l})=g_l(z_{q_{l-1}+1},...,z_{q_l}),~2\leq l\leq t. \\
\end{aligned}
\right.
\end{equation}
Plus the equation $\alpha_1+\alpha_2+\cdots+\alpha_n=g-2k_0\pi i$ in (\ref{2.4}), one gets
\begin{equation}\label{2.8111}
g_1(z_1,...,z_{p_1})+\sum_{l=2}^{t} g_l(z_{q_{l-1}+1},...,z_{q_l})=g-2k_0\pi i.
\end{equation}

Combining (\ref{2.8}), (\ref{2.8111}) and the above argument yields that $u(z_1,...,z_{n})=\psi_1(z_1,...,z_{p_1})+\sum_{l=2}^{t} \psi_l(z_{q_{l-1}+1},...,z_{q_l})$, where $\psi_l$ is an entire function such that
$$
(\psi_1)_{z_1}(\psi_1)_{z_2}\cdots (\psi_1)_{z_{p_1}}=e^{g_1},~~(\psi_l)_{z_{q_{l-1}+1}}(\psi_l)_{z_{q_{l-1}+2}} \cdots (\psi_l)_{z_{q_l}}=e^{g_l},~~2\leq l\leq t.
$$

We observe that the original equation (\ref{2.1}) decomposes into several equations as above with fewer variables. And these equations satisfy the conditions of \textbf{Subcase 1.1}. The conclusion of \textbf{Subcase 1.1} yields
$$
\psi_1=f_1(A_{11}z_1+...+A_{1p_1}z_{p_1}),~~\psi_l=f_l(A_{l1}z_{q_{l-1}+1}+...+A_{lp_2}z_{q_l}),~~ 2\leq l\leq t,
$$
where $f_l$ is entire function satisfying $f_l'=e^{\frac{g_l}{p_l}}$ since $\sharp I_l=p_l$. This is still the conclusion (b).\\

\textbf{Case 2.} $\sharp J>0$. \\

Observe that $J=\{\tau_1,...,\tau_\nu\}$. For each $\tau_i\in J$, then $u_{\tau_i j}\equiv0$ for any $j\in \{1,2,...,n\}\backslash \{\tau_i\}$, which implies that $\frac{\partial \alpha_{\tau_i}}{\partial z_j}\equiv 0$ and $\alpha_{\tau_i}=\alpha_{\tau_i}(z_{\tau_i})$ depends only on $z_{\tau_i}$. The fact $u_{j \tau_i}=u_{\tau_i j}$ yields that $\frac{\partial \alpha_j}{\partial z_{\tau_i}}\equiv 0$ for any $j\in \{1,2,...,n\}\backslash \{\tau_i\}$. Further, $\frac{\partial \alpha_j}{\partial z_{\tau_i}}\equiv 0$ for any $j\in I\backslash J= \chi$ and $\tau_i\in J$. Below, we again consider two subcases.\\

\textbf{Subcase 2.1.} $\sharp J=n$.\\

Then $\chi=\emptyset$ and
\begin{equation}\label{2.41}
\alpha_i=\alpha_i(z_{i}), ~~(i=1,...,n)
\end{equation}
which plus the equation (\ref{2.4}) yields that $g-2 k_0 \pi i=g_1(z_1)+g_2(z_2)+...+g_n(z_n)$ and
$$
\alpha_i=g_i(z_i),~u_{z_i}=e^{\alpha_i}=e^{g_i}.
$$
A routing calculation yields
$$
u=\sum_{i=1}^n\int e^{g_1(z_i)}dz_i,
$$
which is (a). \\

\textbf{Subcase 2.2.} $\sharp J<n$.\\

Then, the fact $\frac{\partial \alpha_j}{\partial z_{\tau_i}}\equiv 0$ for any $j\in I\backslash J=\chi$ yields that

\begin{equation}\label{2.40001}
\alpha_{\tau_{j}}=\alpha_j(z_{\tau_{\nu+1}},z_{\tau_{\nu+2}},...,z_{\tau_n}) ~~\hbox{for} ~~j\in \{\nu+1,...,n\},
\end{equation}
which plus the equation (\ref{2.4}) yields that $g-2 k_0 \pi i=\sum _{i=1}^\nu g_i(z_{\tau_i})+h(z_{\tau_{\nu+1}},z_{\tau_{\nu+2}},...,z_{\tau_n})$ and
$$
\alpha_{\tau_i}=g_i(z_{\tau_i})~~(i=1,...,\nu),~~\sum_{j=\nu+1}^n \alpha_{\tau_j}=h(z_{\tau_{\nu+1}},z_{\tau_{\nu+2}},...,z_{\tau_n}).
$$
In addition,
$$
u=\sum_{i=1}^{\nu} \int e^{g_{\tau_i}(z_{\tau_i})}dz_{\tau_i}+v(z_{\tau_{\nu+1}},z_{\tau_{\nu+2}},...,z_{\tau_n}),
$$
where $v=v(z_{\tau_{\nu+1}},z_{\tau_{\nu+2}},...,z_{\tau_n})$ is an entire function of $n-t$ complex variables satisfying
\begin{equation}\label{2.00005}
v_{z_{\tau_{\nu+1}}}v_{z_{\tau_{\nu+2}}}...v_{z_{\tau_{n}}}=e^{h}.
\end{equation}
Obviously, $v$ and the equation (\ref{2.00005}) satisfy the condition of \textbf{Case 1}. Subsequently, by the conclusion of \textbf{Case 1}, we derive the desired conclusion (c).\\

This finishes the proof of Theorem 2.\\


\begin{thebibliography}{99}

\bibitem{Calin} O. Calin and D. C. Chang, Geometric mechanics on Riemannian manifolds, Applications to partial
differential equations, Birkhauser, Boston, 2005.

\bibitem{Chen} W. Chen and Q. Han, On entire solutions to eikonal-type equations, J. Math. Anal. Appl. 506(2022), 124704.

\bibitem{Courant} R. Courant and D. Hilbert, Methods of Mathematical Physics, II. Partial differential equations, Interscience, New York, 1962.

\bibitem{Garabedian}P.R. Garabedian, Partial differential equations, John Wiley, New York, 1964.

\bibitem{Hem}J.E. Hemmati, Entire solutions of first-order nonlinear partial differential equations, Proc. Am. Math. Soc. 125(1997), 1483-1485.

\bibitem{Hu} P.C. Hu, P. Li and C.C. Yang, Unicity of meromorphic mappings, Advances in Complex Analysis and Its Applications, vol.
1, Kluwer Academic Publishers, Dordrecht, Boston, London, 2003.

\bibitem{Khavinson} D. Khavinson, A note on entire solutions of the eiconal equation, Amer. Math. Monthly. 102(1995), 159-161.

\bibitem{Li1} B.Q. Li, Entire solutions of certain partial differential equations and factorization of partial derivatives, Trans. Amer. Math. Soc. 357(2005), 3169-3177.

\bibitem{Li2} B.Q. Li, On certain functional and partial differential equations, Forum Math. 17(2005), 77-86.

\bibitem{Li3} B.Q. Li, Entire solutions of $(u_{z_1})^m+(u_{z_2})^n=e^g$, Nagoya Math. J. 178(2005), 151-162.

\bibitem{Li4} B.Q. Li, On meromorphic solutions of $f^2 + g^2 = 1$, Math. Z. 258(2008), 763-771.

\bibitem{Lu} F. L\"{u}, Entire solution of a variation of the eikonal equation and related PDEs, P. Edinburgh Math. Soc. 63(2020), 1-12.

\bibitem{Ma} F. L\"{u} and Z.M. Ma, Entire solutions of product type nonlinear partial differential equations in $\mathbb{C}^n$, Glasg. Math. J. Published online 2025:1-6. doi:10.1017/
S0017089525100657.

\bibitem{St} W. Stoll, Introduction to value distribution theory of meromorphic maps, In: Complex Analysis. Lecture Notes in Mathematics, vol 950. Springer, Berlin, Heidelberg. (1982), 210-359.

\bibitem{Saleeby} E.G. Saleeby, Entire and meromorphic solutions of Fermat type partial differential equations, Analysis. 19(1999), 369-376.

\bibitem{Saleeby1} E.G. Saleeby, On entire and meromorphic solutions of $\lambda u^k+\sum_{i=1}^n u^m_{z_i}=1$, Complex Var. Theory Appl. 49(2004) 101-107.

\bibitem{Vi} A. Vitter, The lemma of the logarithmic derivative in seveal complex variables, Duke Math. J. 44(1977), 89-104.

\bibitem{Xu} H.Y. Xu, K. Liu and Z.X. Xuan, Results on solutions of several product type nonlinear partial differential equations in $\mathbb{C}^3$, J. Math. Anal. Appl. 543(2025), ID 128885, 21 p.

\end{thebibliography}
\end{document}